\documentclass[journal,usenames,dvipsnames]{IEEEtran}

\usepackage{myStyleIEEE}
\usepackage{bm}
\usepackage{caption}
\usepackage{smartdiagram}
\usesmartdiagramlibrary{additions}
\begin{document}

\title{Active Distribution Grids Providing Voltage Support: The Swiss Case}

\author{
\IEEEauthorblockN{Stavros~Karagiannopoulos,~\IEEEmembership{Member,~IEEE,} Costas Mylonas,
  Petros~Aristidou,~\IEEEmembership{Member,~IEEE,}\\
  and~Gabriela~Hug,~\IEEEmembership{Senior Member,~IEEE}}%
  \thanks{S. Karagiannopoulos, C. Mylonas and G. Hug are with the Power Systems Laboratory, ETH Zurich, 8092 Zurich, Switzerland. Email: \{karagiannopoulos $|$comylona$|$ hug\}@eeh.ee.ethz.ch.}%
  \thanks{P. Aristidou is with the Department of Electrical Engineering, Cyprus University of Technology, 3036 Limassol, Cyprus. Email: petros.aristidou@cut.ac.cy}%
}

\maketitle

\IEEEpeerreviewmaketitle
\begin{abstract}
The increasing installation of controllable Distributed Energy Resources (DERs) in Distribution Networks (DNs) opens up new opportunities for the provision of ancillary services to the Transmission Network (TN) level. As the penetration of smart meter devices and communication infrastructure in DNs increases, they become more observable and controllable with centralized optimization-based control schemes becoming efficient and practical. In this paper, we propose a centralized tractable Optimal Power Flow (OPF)-based control scheme that optimizes the real-time operation of active DNs, while also considering the provision of voltage support as an ancillary service to the TN. We embed in the form of various constraints the current voltage support requirements of Switzerland and investigate the potential benefit of ancillary services provision assuming different operational modes of the DER inverters. We demonstrate the performance of the proposed scheme using a combined HV-MV-LV test system.
\end{abstract}


\section{Introduction}
Over the last decades, power systems have been facing notable developments in both the transmission and distribution systems. The need to cope with the climate change challenge has given a strong impulse to rely on energy produced by renewable energy resources, phasing out large plants based on fossil fuels connected to the Transmission Networks (TNs)~\cite{gielen2019role}. At the same time, the role of active Distribution Networks (DNs) is greatly upgraded due to the numerous installations of Distributed Energy Resources (DERs) on Medium Voltage (MV) and Low Voltage (LV) levels, as well as due to the vast installation of smart metering devices with control capabilities~\cite{sun2015comprehensive}. \SK{Over the last years, advanced metering infrastructure based on smart meters has gained a lot of attention worldwide. According to~\cite{eia}, $52\%$ of the 150 million electricity consumers in the United States of America have such metering infrastructure and the situation is similar in Europe. In Switzerland, $80\%$ of all electricity meters will have to be replaced with smart meters, e.g.~\cite{landis}, by 2027~\cite{SFOE}. In our work, we consider smart meters that enable two-way communication between the meter and the central system, and therefore can receive remote control signals.}

Distributed Generators (DGs), such as Photovoltaic (PV) units and Wind Turbines (WTs), in combination with other DERs, e.g. electric vehicles, Battery Energy Storage Systems (BESSs) and Controllable Loads (CLs) have changed the paradigm of treating the DNs as sinks of power and offer control capabilities to provide ancillary services and support the bulk transmission system~\cite{PES-TR22}. Thus, to ensure secure system operation coordination and a closer collaboration between Transmission System Operators (TSOs) and Distribution System Operators (DSOs) are necessary.

The efficiency and capability of the DNs to assist TSOs in coping with their challenges depend on technical aspects, such as the available communication and monitoring infrastructure where centralized, distributed and decentralized or local schemes are conceivable, as well as regulatory aspects, e.g. the actual grid codes and standards that describe the allowed ancillary service provision framework. 

The provision of ancillary services by active DNs has started to gain a lot of attention lately~\cite{karagiannopoulos2019active,saint2016active,valverde2019coordination,zerva2014contribution,sun2019review,knezovic2015distribution,valverde2013model}. The applications range from frequency regulation~\cite{karagiannopoulos2019active} and congestion management~\cite{knezovic2015distribution} to voltage support~\cite{valverde2013model,zerva2014contribution, valverde2019coordination,sun2019review} which is the focus of this paper. \SK{Existing works on voltage control differ in terms of the available information and communication infrastructure, as well as the existing control architecture. First, purely \emph{local} control schemes, e.g.~\cite{valverde2019coordination} are the cheapest and most scalable alternative in grids with numerous DERs and without communication capabilities. Then, \emph{distributed} approaches, e.g.~\cite{cavraro2016value,robbins2012two,arnold2017model}, allow neighbors to communicate values in order to coordinate reactive power exchanges, and finally, \emph{centralized} hierarchical, e.g.~\cite{abessi2015centralized}, and optimization-based methods, e.g.~\cite{valverde2013model,bidgoli2016receding,saint2016active}, control the response of all units to optimize the overall system's objectives.} In~\cite{valverde2013model}, a centralized approach based on model predictive control is used to track voltages in a MV grids. The voltage sensitivities are used instead of the power flow equations and only the MV level is considered. Another centralized approach is presented in~\cite{saint2016active}, where an OPF-based dual-horizon rolling scheduling model calculates first the optimal schedule for the power exchange between TN and DN and minimizes the real-time deviations in the operation stage. However, emphasis is put on the active power exchange, the DGs are modelled as non-dispatchable active power injections without considering the reactive power capabilities by the inverter control, and only the MV grid is modeled. Reference~\cite{zerva2014contribution} investigates the active and passive participation in the voltage support scheme of Switzerland~\cite{SwissGrid_VoltageSupport} from the TSO perspective. The same concept is examined in~\cite{valverde2019coordination} which presents a model-free control scheme for DERs located in LV DNs to provide voltage support to TNs. The default local control scheme of the DERs is altered upon request to support the TN voltage and the case studies consider the impact of both MV and LV DERs. However, by design, only the DERs which do not face any power quality issue will contribute to voltage support, e.g. the ones close to the substation where voltage is regulated. Another model-less approach is presented in~\cite{arnold2017model} and uses the theory of extremum seeking to optimally track a voltage profile at the TN and DN interconnection. Nevertheless, tuning and implementation challenges impose constraints on this promising approach. Finally, a detailed review of voltage control schemes in TNs and DNs is presented in~\cite{sun2019review} discussing open topics and challenges in the TN-DN coordination potential. \SK{In contrast to most references which focus on the modeling solely of the MV~\cite{valverde2013model,bidgoli2016receding} or LV~\cite{valverde2019coordination,cavraro2016value} grid, we investigate simultaneously the MV and LV grids in order to consider the contribution of the numerous DERs in LV grids. Additionally, we incorporate the BFS power flow equations and we do not rely on voltage sensitivities~\cite{valverde2013model}.}

This paper presents a centralized tractable OPF-based control scheme that optimizes the real-time operation of active DNs, considering also the provision of voltage support as an ancillary service to the TN. The proposed tool is capable of accommodating any inverter-based DER in MV and LV, such as WTs, solar parks, rooftop PV units, BESS, and CLs.
\SK{Traditionally, passive distribution networks have relied on conventional network elements such as capacitor banks, static var compensators and on-load tap changers, or they impose strict limits and bounds for demand and generation units. Such technologies can be easily incorporated in our formulation and models, e.g. a static var compensator can be modeled by setting the active power exchange to zero and the reactive power constrained by the technical minimum and maximum values. Other DER technologies, such as synchronous generators and doubly-fed induction generators for wind turbines, can also be easily included in the formulation. 
However, in this paper, we focus on demonstrating the ability of DERs, which are expected to play a key role in active distribution grids, to assist the DSOs not just in active power balancing but also in voltage control.} Tractability is achieved by formulating the controller as a Backward/Forward Sweep (BFS) OPF, extending our previous works~\cite{StavrosIREP,StavrosPSCC18,karagiannopoulos2019data}. In contrast to previous work, we investigate the technical potential of the inverters under various existing grid codes and standards, and we highlight the need of harmonization and modernization of new guidelines. Furthermore, we consider the HV grid as a Th\'{e}venin equivalent allowing us to examine the impact that active control in MV and LV grids can have on the TN, highlighting the importance of combined analysis and TN-DN coordination and joint planning. \SK{In the future, large shares of operational flexibility will be coming from DERs that are connected in LV grids (electric vehicles, batteries, rooftop PV panels, etc.). Overvoltages due to high PV injections, congestion in the cables due to synchronized charging of electric vehicles or batteries, or even unbalances due to different behavior among the phases are some problems that cannot be mitigated by OPF-based control without detailed LV grid modeling.}

\SK{The contributions of this paper can be summarized as follows:
\begin{itemize}
    \item Detailed modeling of both the MV and LV grids in order to guarantee safe grid operation on all voltage levels. The HV network is represented by a Th\'{e}venin equivalent according to current practice.
    \item Focus on inverter-based DERs to investigate the coordinated potential of DERs that are expected to play a key role in active distribution grids.
    \item Investigation of different feasible areas in terms of the inverter P-Q characteristic in order to highlight the full potential of reactive power control and the need for new grid codes.
    \item Formulation of a computationally tractable BFS-OPF considering the  provision of voltage support ancillary service.
    \item Mathematical formulation of the voltage support scheme in Switzerland, and monetary investigation for both participation types through a case study.
\end{itemize}}

The remainder of the paper is organized as follows: In Section~\ref{VSinCH}, we present the reasoning and the requirements to participate in the voltage control scheme in Switzerland, either having a passive or an active role. Then, in Section~\ref{detOPF}, we give the formulation of the centralized BFS-OPF with emphasis on the constraints for the voltage support case and the different inverter operating regions. In Section~\ref{case}, we introduce the case study and simulation results that show the performance of the optimized controllers. Finally, we draw conclusions in Section~\ref{Conclusion}.

\section{Voltage Support Schemes in Switzerland}\label{VSinCH}
Traditionally, DNs were treated as sinks of power with limited demand uncertainty. Modern DNs, however, host various types of DERs with both active and reactive power control capabilities, enabling more efficient and complex TSO-DSO interactions. In this work, we investigate the provision of voltage support from active DNs in Switzerland. In this section, we introduce the two available support schemes, namely the \emph{passive} and the \emph{active} types.

\subsection{Passive Scheme}
\begin{figure}[]
    \centering
    \begin{minipage}{.25\textwidth}
        \centering
        \includegraphics[width=1\textwidth]{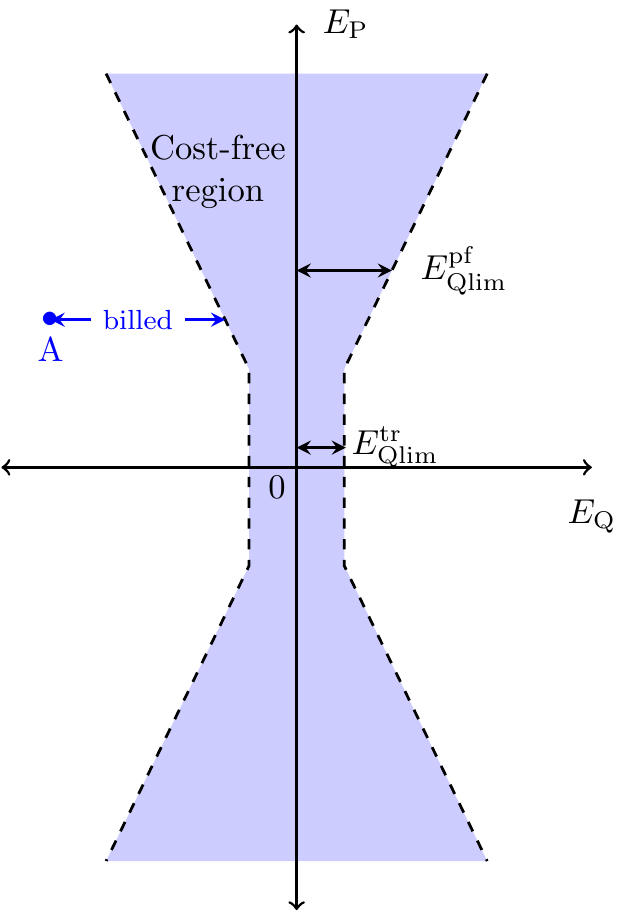} 
        \captionsetup{labelformat=empty}
        \caption*{(a)}
        \label{fig:pas}
    \end{minipage}%
    \begin{minipage}{.25\textwidth}
        \centering
        \includegraphics[width=1\textwidth]{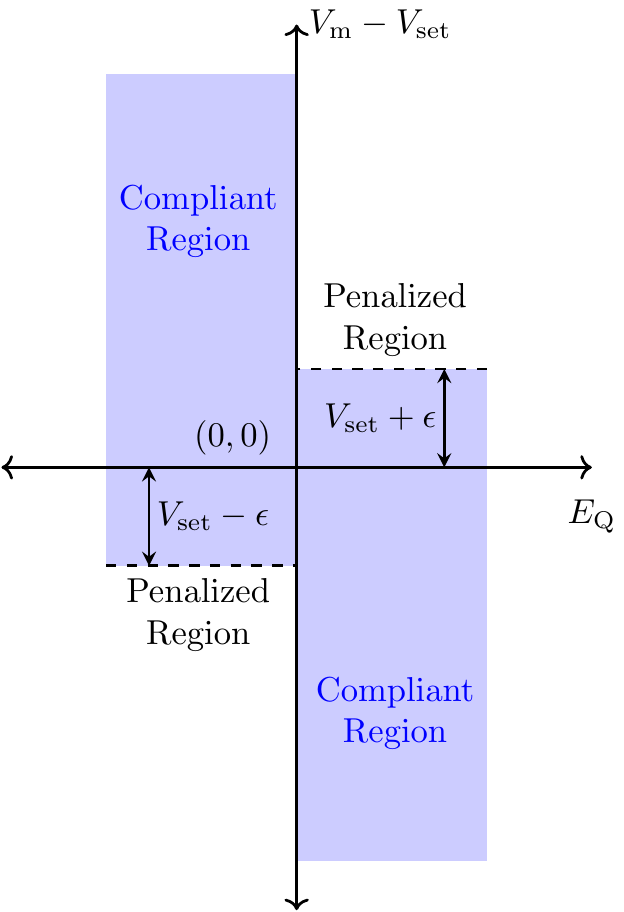} 
        \captionsetup{labelformat=empty}
        \caption*{(b)}
        \label{fig:act}
    \end{minipage}
    \vspace{+0.1cm}
    \caption{Passive and active participation in the voltage support scheme of the Swiss TN~\cite{SwissGrid_VoltageSupport}.
	(a) cost-free region of the passive participation based on the 15-min energy exchange between the DN and TN; (b) compliant and non-compliant regions of the active participation.}
    \label{fig:passive_scheme}
\end{figure}

The main scope of the passive voltage support scheme is to impose a limit on the impact of the DN operation on the TN voltage at the transmission-distribution connection point, i.e. at the primary side of the HV-MV transformer. This is achieved by obliging the DN to operate with a minimum power factor that depends on the active power exchange. The symmetrical cost-free region in the active and reactive energy exchange plane of this scheme is shown in Fig.~\ref{fig:passive_scheme}(a).    

The values for the active and reactive energy exchange with the TN, namely $E_{\textrm{P}}$ and $E_{\textrm{Q}}$, are measured over a 15-minute window and the DN is penalized based on the deviation from the cost-free region. \SK{The boundaries of this region are defined by the predetermined minimum power factor $\cos(\phi_{\textrm{min}})$ which defines the cost-free reactive energy exchange as 
\begin{equation} \label{eq:W_Q_lim_LF}
|E_{\textrm{Qlim}}^{\textrm{pf}}| \leq \tan(\phi_{\textrm{min}}) \cdot  E_{\textrm{P}}.  
\end{equation}
\noindent where $E_{\textrm{Qlim}}^{\textrm{pf}}$ is the reactive energy exchange limit due to the minimum power factor requirement shown in Fig.\ref{fig:passive_scheme}(a).}

In the past, DSOs preferred disconnecting the transformer at very lightly loaded conditions, i.e. when limited active energy exchange with the TN was needed. The loads were supplied by other transformers, e.g. in a ring-shaped MV configuration, and the DSOs avoided penalization in presence of reactive energy exchanges. However, such strategical operation influenced the security of supply and therefore, the boundary is enlarged for small values of active power needs, given by 
\begin{equation} \label{eq:W_Q_lim_TRAFO}
|E_{\textrm{Qlim}}^{\textrm{tr}}| \leq \frac{u_{k}}{100} \cdot S_{N} \cdot \Delta t_V,
\end{equation}
\noindent where $u_{k}$ is the transformer's short-circuit voltage (\%), $S_{\textrm{N}}$ its nominal apparent power and $\Delta t_V$ the 15-minute measurement period.

Overall, the final cost-free reactive energy exchange is given by
\begin{equation}
E_{\textrm{Qlim}} = max\{ E_{\textrm{Qlim}}^{\textrm{pf}}, E_{\textrm{Qlim}}^{\textrm{tr}}\}.,
\end{equation}

The cost of the net energy values which are measured every 15 minutes for the passive participation case ($C_{\textrm{P}}$) is given by

\begin{equation} \label{eq:W_Q_ver}
\centering
  C_{\textrm{P}}=\left\{
  \begin{array}{@{}ll@{}}
    c_{\textrm{p}} ( |E_{\textrm{Q}}| - E_{\textrm{Qlim}}), & \text{if}\ |E_{\textrm{Q}}| > E_{\textrm{Qlim}}  \\
    \qquad 0, & \text{otherwise}
  \end{array}\right.,
\end{equation} 
\noindent where $c_{\textrm{p}}$ is the reactive energy tariff in $\frac{CHF}{Mvarh}$ for the penalized region ~\cite{SwissGrid_VoltageSupport_Tariff}. \SK{For 2021, this cost is set to $C_{\textrm{P}}=0.0138\frac{CHF}{Mvarh} \approx 0.0138 \frac{\$}{Mvarh}$.} 



\subsection{Active Scheme} \label{activeVscheme}
In the active voltage support scheme, the participating DSOs have to follow a time-varying voltage reference $V_{set}$ provided by the TSO for the substation bus connecting the DN to the TN. These setpoints are calculated off-line, based on the Day-ahead Reactive Planning (DARP) problem, i.e. an OPF problem using forecast data with an hourly time resolution~\cite{zerva2014contribution}. The single-period OPF model is solved for each hour of the next day using data and models available to the TSO, and defines the setpoints for the tap positions of the transformers, and the hourly voltage setpoints that the active participants need to follow the next day.


Figure~\ref{fig:passive_scheme}(b) depicts the compliant regions for the active participation in order to track the voltage setpoint $V_{set}$. \SK{The DSO that participates in the voltage support ancillary service is compliant when its reactive energy exchange is assisting the TSO to track $V_{set}$.} \SK{That is, providing reactive power ($E_Q>0$) when the measured voltage is below the reference value ($V_{\textrm{m}}<V_{\textrm{set}}$), and consuming reactive power when the measured voltage is above the reference value ($V_{\textrm{m}}>V_{\textrm{set}}$). Please note that reference~\cite{SwissGrid_VoltageSupport} uses the opposite sign convention following the consumer's reference convention, i.e. positive reactive power corresponds to an inductive behavior reducing the local voltage magnitude value.} A small deviation tolerance with respect to the voltage setpoint value is allowed, e.g. for connection points to the $380~kV$ grid the measurement tolerance is $\Delta V_{\textrm{m}} = 3~kV$ and for the $220~kV$ grid $\Delta V_{\textrm{m}} = 2~kV$~\cite{swissgrid}. The compliant ($A_1$, $A_2$) and non-compliant ($A_3$, $A_4$) regions are described as follows:
\SK{\begin{subequations} 
\begin{align}
& A_{1} = \{E_{\textrm{Q}} \leq 0, V_{\textrm{m}}-V_{\textrm{set}} \geq -\epsilon \}, \\
& A_{2} = \{E_{\textrm{Q}} \geq 0, V_{\textrm{m}}-V_{\textrm{set}} \leq \epsilon \}, \\
& A_{3} = \{E_{\textrm{Q}} \leq 0, V_{\textrm{m}}-V_{\textrm{set}} \leq -\epsilon \},
\\
& A_{4} = \{E_{\textrm{Q}} \geq 0, V_{\textrm{m}}-V_{\textrm{set}} \geq \epsilon  \}. 
\end{align} 
\end{subequations}}
The financial compensation for the reactive energy exchanged depends on the 15-minute-interval compliance score. If the measured voltage lies within the compliant region, the \SK{DSO} receives a remuneration and conversely if it lies in the non-compliant region a penalty is applied. The final cost for the DSO is proportional to the amount of reactive energy exchanged, given by

\begin{equation} \label{eq:VB_Q_active}
  C_{\textrm{A}}=\left\{
  \begin{array}{@{}ll@{}}
    -c_{\textrm{a,c}} \cdot \lvert E_{\textrm{Q}} \rvert , & \text{if}\ E_{\textrm{Q}} \in A_{1} \cup A_{2}  \\
    +c_{\textrm{a,n}} \cdot \lvert E_{\textrm{Q}} \rvert , & \text{if}\ E_{\textrm{Q}} \in A_{3} \cup A_{4} \
  \end{array}\right.,
\end{equation} 
where $c_{\textrm{a,c}}$ is the conforming revenue tariff and $c_{\textrm{a,n}}$ the non-conforming penalty for the 15-minute remuneration period~\cite{SwissGrid_VoltageSupport_Tariff}. If the \SK{DSO} is compliant at least 80\% of the time over a month, it is remunerated. However, monthly compliance below 70\% for two consecutive months automatically results in a switch to a passive participation for the DN in question. \SK{For 2021, these costs are set to $c_{\textrm{a,c}}=0.003\frac{CHF}{Mvarh} \approx 0.003\frac{\$}{Mvarh}$ and $c_{\textrm{a,n}}=0.0138\frac{CHF}{Mvarh} \approx 0.0138\frac{\$}{Mvarh}$.} 

\section{Centralized Optimal Power Flow Formulation}\label{OPFall}
As mentioned previously, in this work, we assume that the DSO has the ability to dispatch DERs in a centralized way through an existing communication system. In this section, we present the centralized OPF-based control scheme that drives the DER setpoints to satisfy the power quality and security constraints, as well as provide voltage support to the TN. Finally, the OPF formulation allows us to compare the performance and quantify the benefits of using different operational types of inverter-based DERs.
\subsection{Centralized OPF}\label{detOPF}
\subsubsection{Objective function}\label{objfun}
The objective function minimizes the cost of DER control and the network losses, over all of the network nodes ($N_b$) and branches ($N_{br}$) for the entire control time horizon ($N_{hor}$). Moreover, it optimizes the costs from providing voltage support to the TN, either as a passive or active participant. This is described by: 
\begin{align}
    \min_{ \bm{u} } & \sum \limits_{t=1}^{N_\textrm{hor}} \biggl\{  \sum \limits_{j=1}^{N_\textrm{b}}  \biggl(C_\textrm{curt} \hspace{-0.05cm} \cdot \hspace{-0.05cm} P_\textrm{curt,j,t} \hspace{-0.05cm} + \hspace{-0.05cm} C_\textrm{Q} \hspace{-0.05cm} \cdot \hspace{-0.05cm} Q_\textrm{ctrl,j,t} \biggr) \cdot \Delta T + \nonumber  \\  
    & \sum \limits_{i=1}^{N_\textrm{br}} C_\textrm{curt} \hspace{-0.05cm} \cdot \hspace{-0.05cm} P_\textrm{loss,i,t} \cdot \Delta T + \biggl(   C_\textrm{P,t} +  C_\textrm{A,t} \biggl) \biggl\} \nonumber \\
    &+ C_\textrm{H} \cdot \biggl( ||\eta_\textrm{V}||_{\infty} + ||\eta_{\textrm{I}}||_{\infty}  \biggr), \label{eq:objfun1} 
\end{align}
\noindent where $\bm{u}$ is the vector of the available control measures and $\Delta T$ is the length of each time period (in this work 15 minutes). The curtailed power of the DGs connected at node $j$ and time interval $t$ is given by $P_\textrm{curt,j,t} = P_\textrm{g,j,t}^{\textrm{max}} - P_\textrm{g,j,t}^{\textrm{ }}$, where $P_\textrm{g,j,t}^{\textrm{max}}$ is the maximum available active power and $P_{\textrm{g,j,t}}^{\textrm{ }}$ the active power injection of the DGs. \SK{Thus, we account for consumers' compensation in case their production needs to be curtailed with the cost $C_\textrm{curt}$.} The use of reactive power support $Q_\textrm{ctrl,j,t}= | Q_\textrm{g,j,t}^{\textrm{ }}|$ for each DG connected to node $j$ and time interval $t$ is also minimized where $Q_{\textrm{g,j,t}}^{\textrm{ }}$ is the DG reactive power injection or absorption. The coefficients $C_\textrm{curt}$ and $C_\textrm{Q}$ represent the cost of curtailing active power and providing reactive power support (DG opportunity cost or contractual agreement). The assumption that $C_\textrm{Q} \ll C_\textrm{curt}$ is made, which prioritizes the use of reactive power control over active power curtailment. 

\SK{This paper takes the perspective of the DSO which can use part of the operational flexibility provided by DERs. Thus, it is assumed that the DER owners and the DSO have agreed to such an intervention in exchange for some economical benefit or as part of their licensing agreement. The units that do not participate in this scheme, are treated as inelastic from the DSO, i.e. the DSO has no control over their behavior apart from making sure they operate according to the existing grid code without violating the prescribed power quality constraints.}

\SK{To provide non-discriminatory access to all customers and treat them in a fair way, the same compensation remuneration or cost penalization schemes are used for all customers irrespective of their location in the network or their actual different production costs.}

The \SK{per unit losses assuming balanced loading} in each branch $i$ at time $t$ are calculated by $P_\textrm{loss,i,t} = |I_{\textrm{br,i,t}}| ^{2} \hspace{-0.05cm} \cdot \hspace{-0.05cm} R_{\textrm{br,i}}$, where $|I_\textrm{br,i,t}|$ is the magnitude of the current flow and $R_{\textrm{br,i}}$ the resistance of the branch. \SK{We include the loss minimization aspect, since it is recognized as an important opportunity for active DSOs which will be able to utilize the operational DG flexibility~\cite{farrokhseresht2017minimization,EDSO}.}

The costs of the passive participation case ($C_\textrm{P,t}$) and the costs or financial remuneration in case of active participation in the voltage support schemes ($C_\textrm{A,t}$) are considered for each time interval $t$ in the objective function.


Finally, $C_\textrm{H}$ is a large penalty linked with security constraint and power quality constraint violations. It is used in conjunction with the variables $(\eta_\textrm{V},\eta_{\textrm{I}})$ to relax respectively the voltage and thermal constraints and avoid deriving infeasible solutions. When one of these limits is violated, the output of the overall objective function is dominated by this term and might lose a real monetary meaning (unless the cost of violating the security and power quality constraints is quantified and monetized by the DSO).
\subsubsection{Power balance constraints}\label{powerbalance}
The power injections at every node $j$ and time step $t$ are given by
\begin{subequations} \label{eq:node_balance}
\begin{align}
	P_{\textrm{inj,j,t}}^{\textrm{ }}&=P_{\textrm{g,j,t}}^{\textrm{ }} - P_{\textrm{lflex,j,t}}^{\textrm{ }} - (P_{\textrm{B,j,t}}^{\textrm{ch}} - P_{\textrm{B,j,t}}^{\textrm{dis}}), \label{eq:node_balance_P}\\
    Q_{\textrm{inj,j,t}}^{\textrm{ }}&=Q_{\textrm{g,j,t}}^{\textrm{ }} - Q_{\textrm{lflex,j,t}}^{\textrm{ }}  + Q_{\textrm{B,j,t}},  \label{eq:node_balance_Q}   
\end{align}
\end{subequations}
where $P_{\textrm{lflex,j,t}}^{\textrm{ }}$ and $Q_{\textrm{lflex,j,t}}^{\textrm{ }}$ are the active and reactive node demands (after control) of constant power type; $Q_{\textrm{B,j,t}}$ the reactive power of the BESS and, $P_{\textrm{B,j,t}}^{\textrm{ch}}$ and $P_{\textrm{B,j,t}}^{\textrm{dis}}$ are respectively the charging and discharging active powers of the BESS.
\subsubsection{Power flow constraints}\label{powerflow}
The non-linear AC power-flow equations that model the DN network make solving the OPF problem computationally challenging. Since the OPF will be used to process several scenarios in a multi-period framework, it is necessary to use some approximations to increase its computational performance. For this reason, the iterative BFS power flow~\cite{Teng2003} method is used in this work, extending the formulation presented by the authors in \cite{stavrosPowertech,StavrosIREP,StavrosPSCC18,karagiannopoulos2019data} for a single and three-phase systems.
Following our previous work, a single iteration of the BFS power-flow method is used to replace the AC power-flow constraints in the OPF formulation. This is written as ($j=1,\ldots N_{b}$):
\begin{gather}
    I_\textrm{inj,j,t}= \left (\frac{(P_{\textrm{inj,j,t}}^{\textrm{ }} + jQ_{\textrm{inj,j,t}}^{\textrm{ }})^{*}}{\bar{V}_{\textrm{j,t}}^{*}}\right), \nonumber  \\
    I_\textrm{br,t}=BIBC \cdot I_\textrm{inj,t}, \nonumber \\
    \Delta{V}_{\textrm{t}}=BCBV \cdot I_\textrm{br,t}, \nonumber  \\
    V_{{\textrm{t}}}=V_\textrm{slack} - \Delta V_\textrm{tap} \cdot \rho_{\textrm{t}} + \Delta{V}_{\textrm{t}}, \nonumber \\
    \rho_{min} \leq \rho_{\textrm{t}} \leq \rho_{max}, \label{eq:OLTC3}  
\end{gather}
\noindent where $\bar{V}_{\textrm{j,t}}^{*}$ is the voltage at node $j$ at time $t$,$~^*$~indicates the complex conjugate and the bar indicates that the value from the previous BFS iteration is used; $I_\textrm{inj,t}^{\textrm{}}$ and $I_\textrm{br,t}^{\textrm{}}$ are respectively the bus injection and branch flow currents; and, $BIBC$ (Bus Injection to Branch Current) is a matrix with ones and zeros, capturing the topology of the DN; $\Delta{V}^{}_{\textrm{t}}$ contains the voltage drops over all branches and phases; $BCBV$ (Branch Current to Bus Voltage) is a matrix with the complex impedance of the lines as elements; $V_\textrm{slack}$ is the per unit voltage at the slack bus (here assumed to be $1\hspace{-0.1cm}<\hspace{-0.1cm}0^{\circ{}}$ p.u.); $\Delta V_{tap}$ is the voltage magnitude change caused by one tap action of the On-Load Tap Changing (OLTC) transformer and assumed constant for all taps for simplicity; and, $\rho_{\textrm{t}}$ is an integer value defining the OLTC position. The parameters ($\rho_\textrm{min},\rho_\textrm{max}$) are respectively the minimum and maximum tap positions of the OLTC transformer.
This convex formulation provides a good approximation of the nonlinear AC OPF~\cite{Fortenbacher2016a}, is computationally tractable even in a three-phase model~\cite{StavrosPSCC18,karagiannopoulos2019data}, and results in AC feasible solutions which can account for uncertainties, see~\cite{StavrosIREP}. \SK{More specifically,~\cite{Fortenbacher2016a} which is based on a purely linear formulation of the BFS-OPF reports a $2\%$ difference in terms of the objective value compared to solving the exact AC OPF using interior-point methods. This value is calculated for one iteration of the scheme, and the formulation does not consider reactive power flows.}

\begin{figure}[b]
\centering
	\includegraphics[width=0.9\columnwidth]{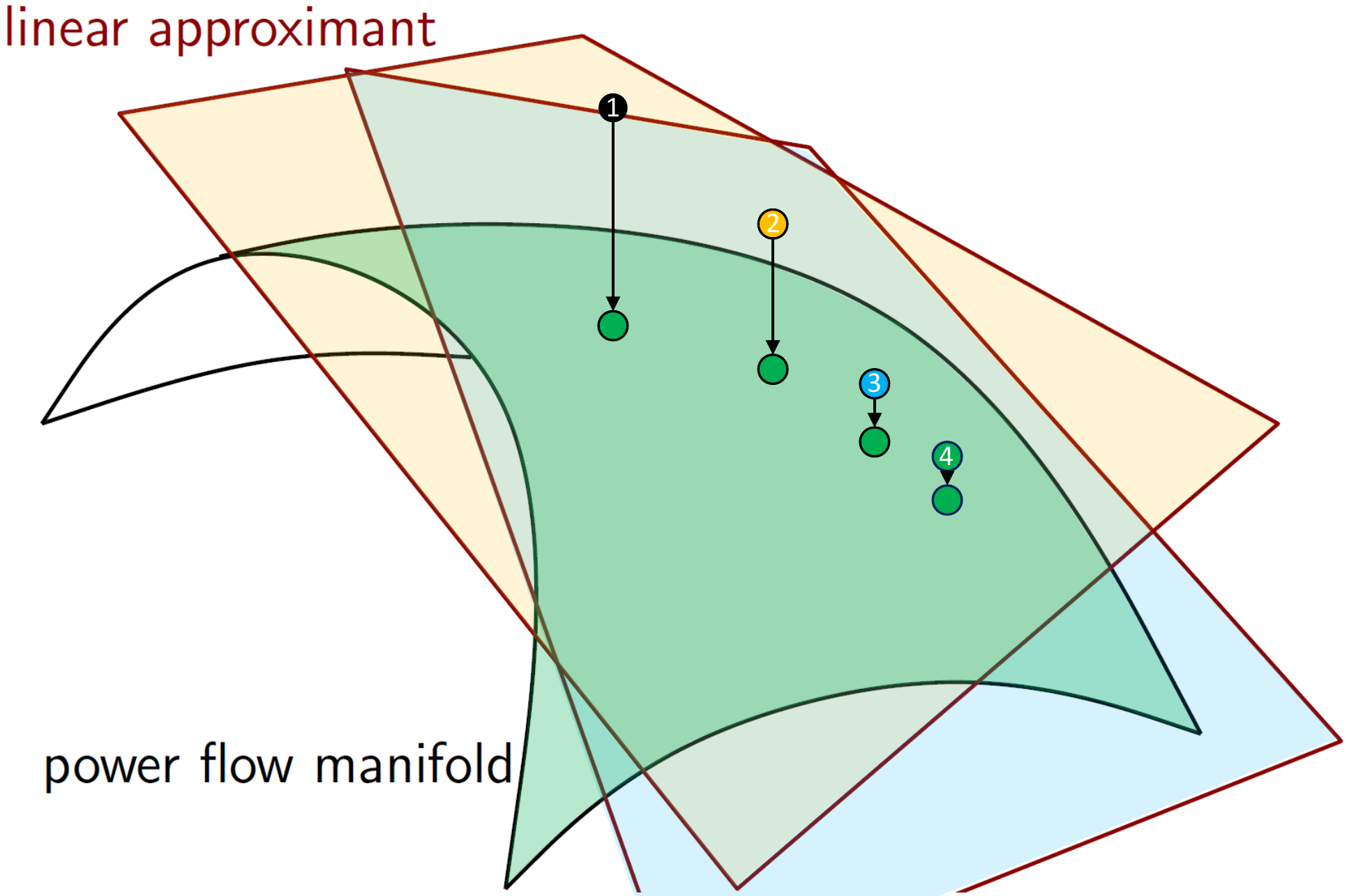}
	\caption{Qualitative illustration of the iterative BFS-OPF scheme.} 
	\label{Ch3manifold}
\end{figure}

In the following, we provide a qualitative illustration of the solution procedure for the BFS-OPF calculations. Inspired by~\cite{Bolognani_manifoldlinearization}, we illustrate in Fig.~\ref{Ch3manifold} the iterative scheme used in the BFS-OPF formulations. The solution of the OPF problem using one BFS iteration does not in general lie on the AC power flow manifold, i.e., does not satisfy the exact power flow equations which is shown in green. Hence, the derived OPF solutions are shown with the operating points denoted by numbers indicating the iteration index. Then, by running an exact PF solution using the optimal setpoints, we project the solution onto the feasible AC manifold and repeat the procedure until convergence. Finally, in this example the algorithm is terminated after 4 iterations, where the derived OPF solution lies very close to the exact AC feasible domain. \SK{In our algorithm, convergence is achieved when the maximum voltage magnitude deviation between the exact power flow solution and the optimization outcome is below a threshold of $10^{-4}$~\cite{StavrosIREP,karagiannopoulos2019data}.} 

\subsubsection{Passive voltage support constraints}\label{Unbalance}
The constraints for the cost-free region of the passive participation in the Swiss voltage scheme are given by
\begin{align}
  - M \cdot Z_\textrm{P,t} + E_{\textrm{Qlim}}  \leq 
|E_{\textrm{Q,t}} | \leq  E_{\textrm{Qlim}}+ M \cdot (1-Z_\textrm{P,t}),
    \label{eq:passive} 
\end{align}
\noindent where $E_{\textrm{Q,t}} =Q_{\textrm{inj,$f_p$,t}}\cdot \Delta t_{V}$ is the reactive power exchange with the TN at time $t$,  $E_{\textrm{Qlim}}=\max \{ E_{\textrm{Qlim}}^{\textrm{pf}}, E_{\textrm{Qlim}}^{\textrm{tr}}\}$, the limit of the cost-free passive participation scheme, $Z_\textrm{P,t}$ the binary variable of the passive role that becomes 1 when the reactive power exchange at time $t$ lies in the cost-free region and $0$ otherwise, and $M$ a sufficiently large constant. \SK{The selection of this constant can have an influence on the final computational time and the ability of a solver to reach an optimal solution. In this paper, setting $M=100$ resulted in a numerically stable solver behavior for all cases examined.} The cost for the passive participation at time $t$ is given by $ C_\textrm{P,t}= (1-Z_\textrm{P,t}) \cdot  (c_{\textrm{p}} ( |E_{\textrm{Q}}| - E_{\textrm{Qlim}}))$.

\subsubsection{Active voltage support constraints}\label{Active}
The constraints for the active voltage support scheme are also modelled using binaries in order to define the compliant and non-compliant regions. In this case, we need one binary to indicate the direction of the cost-free region in terms of reactive power and one to position the measured voltage relative to the reference value considering the tolerance. Thus, the active participation constraints are given by
\begin{gather}
  - M \cdot ( 1 - Z_\textrm{$A_Q$,t})   \leq E_{\textrm{Q,t}}  \leq M \cdot Z_\textrm{$A_Q$,t}, \\
    - M \cdot ( 1 - Z_\textrm{$A_V$,t}) + V_{\textrm{set,t}} - \epsilon  \leq V_{\textrm{m,t}}  \leq M \cdot Z_\textrm{$A_V$,t} + V_{\textrm{set,t}} - \epsilon,
\label{eq:active} 
\end{gather}
\noindent where $V_{\textrm{m,t}}$ is the measured voltage magnitude at the TN-DN interconnection point, $V_{\textrm{set,t}}$ its reference value to be tracked calculated by the TSO, $Z_\textrm{$A_Q$,t}$ is the binary that becomes 1 for positive reactive power exchange with the TN at time $t$ and 0 otherwise, and $Z_\textrm{$A_V$,t}$ the binary that becomes 1 (respectively 0) when the measured voltage is above (respectively below) the reference value taking also the tolerance $\epsilon$ into account.
The compliant regions correspond to  $Z_\textrm{$A_Q$,t}==Z_\textrm{$A_V$,t}$, i.e. when both the binaries are either $1$ or $0$. Therefore,  the revenues and costs operating in the compliant and, respectively non-compliant, area  can be distinguished by these two binaries.  Thus, the costs for the active voltage support participation at time $t$ is given by $C_\textrm{A,t} = \min(Z_\textrm{$A_V$,t},Z_\textrm{$A_Q$,t}) \cdot (-c_{\textrm{a,c}} \cdot | E_{\textrm{Q}}|) + (1- \min(Z_\textrm{$A_V$,t},Z_\textrm{$A_Q$,t}) \cdot (c_{\textrm{a,n}} \cdot |E_{\textrm{Q}}|$).


\subsubsection{Thermal loading and voltage constraints} \label{VIconstraints}
The constraint for the current magnitude for  branch $i$ at time $t$ is given by
\begin{align}
    |I_\textrm{br,i,t}|   \leq I_{\textrm{i,max}} + \eta_{\textrm{I,i,t}} \label{eq:Qk1}, \qquad 
    \eta_{\textrm{I,i,t}}  \geq 0,
\end{align}
\noindent where $I_\textrm{br,i,t}$ is the branch current; $I_{\textrm{i,max}}$ is the maximum thermal limit; and, $\eta_{\textrm{I,i,t}}$ is used to relax the constraint to avoid infeasibility. Similarly, the voltage constraints are given by
\begin{align}
    V_\textrm{min} - \eta_{\textrm{V,j,t}} & \leq | V_{\textrm{j,t}} | \leq V_{\textrm{max}} +\eta_{\textrm{V,j,t}} \label{eq:V1}, \qquad
    \eta_{\textrm{V,j,t}} & \geq 0,
\end{align}
\noindent where $V_{\textrm{min}}$ and $V_{\textrm{max}}$ are the lower and upper acceptable voltage limits, respectively, and $\eta_{\textrm{V,j,t}}$ is used to relax the constraint if necessary.
Unfortunately, \eqref{eq:V1} is non-convex due to the minimum voltage magnitude requirement. In order to avoid the non-convexity, we relax the minimum voltage constraint  (see~\cite{StavrosPSCC18} for more details for the three-phase case), as follows
\begin{align}
    \begin{cases}
               |V_{\textrm{j,t}} | \leq V_{\textrm{max}} + \eta_{\textrm{V,j,t}} \\
               \textrm{Re}\left\{V_{\textrm{j,t}}    \right\} \geq V_\textrm{min} - \eta_{\textrm{V,j,t}} 
            \end{cases}.\label{eq:V4}
\end{align}
\subsubsection{DER constraints}\label{DERineq}
\paragraph{DG limits}
Without loss of generality, we only consider inverter-based DGs such as PVs and WTs. Their limits are thus given by
\begin{subequations}
\label{eq:prod3}
\begin{align}
    P_{\textrm{g,j,t}}^{\textrm{min}} \leq P_{\textrm{g,j,t}}^{\textrm{ }} \leq P_{\textrm{g,j,t}}^{\textrm{max}},\\ Q_{\textrm{g,j,t}}^{\textrm{min}} \leq  Q_{\textrm{g,j,t}}^{\textrm{ }} \leq Q_{\textrm{g,j,t}}^{\textrm{max}},
\label{eq:PV_prod3Qa}
\end{align}
\end{subequations}
\noindent where $P_{\textrm{g,j,t}}^{\textrm{min}}$, $P_{\textrm{g,j,t}}^{\textrm{max}}$, $Q_{\textrm{g,j,t}}^{\textrm{min}}$ and $Q_{\textrm{g,j,t}}^{\textrm{max}}$ are the lower and upper limits for active and reactive DG power at each node $j$,  and time $t$. These limits vary depending on the type of the DG and the control schemes implemented. The focus of this paper is to investigate the benefits of enlarging the allowed region of inverter operational regions, and thus, the next section will provide a detailed modelling of the different alternatives for \eqref{eq:prod3}.
\paragraph{Controllable loads} 
We also consider flexible loads which can shift a fixed amount of energy consumption in time. The behavior of the controllable loads we model is given by
\begin{subequations}\label{eq:CL}
\begin{align}
    P_{\textrm{lflex,j,t}}^{\textrm{ }} &= P_{\textrm{l,j,t}}^{\textrm{ }} + n_{\textrm{j,t}} \cdot P_{\textrm{shift,j}}, \\
   \sum \limits_{t=1}^{N_{hor}} n_{\textrm{j,t}}&=0, \label{eq:CL2}
\end{align}
\end{subequations}
\noindent where $P_{\textrm{lflex,j,t}}^{\textrm{ }}$ is the controlled active power demand at node $j$ and at time $t$, $P_{\textrm{shift,j}}$ is the load that can be shifted (assumed constant) and $n_{\textrm{j,t}} \in \left\{-1,0,1\right\}$ is an integer variable indicating an increase or a decrease of the load when shifted from the initial demand $P_{\textrm{l,j,t}}^{\textrm{ }}$. \SK{It is imposed by~\eqref{eq:CL2} that the final total daily energy demand needs to be maintained.}
\paragraph{Battery Energy Storage Systems}
The constraints related to the BESS are given as
\begin{subequations}
\label{eq:BESS}
\begin{gather}
    SoC_{\textrm{min}}^{\textrm{bat}} \cdot E_{\textrm{cap,j}}^{\textrm{bat}} \leq  E_{\textrm{j,t}}^{\textrm{bat}} \leq SoC_{\textrm{max}}^{\textrm{bat}} \cdot E_{\textrm{cap,j}}^{\textrm{bat}},      \label{eq:BESS_en}\\  
    E_{\textrm{j,1}}^{\textrm{bat}} = E_{\textrm{start}} = \SK{E_{\textrm{j,$N_{hor}$}}^{\textrm{bat}}},    \label{eq:BESS_SOC}\\
    0 \leq P_{\textrm{B,j,t}}^{\textrm{ch}} \leq P_{\textrm{max}}^{\textrm{bat}} ,  \quad  0 \leq P_{\textrm{B,j,t}}^{\textrm{dis}} \leq P_{\textrm{max}}^{\textrm{bat}},  \label{eq:BESS_P} \\
    E_{\textrm{j,t}}^{\textrm{bat}}  = E_{\textrm{j,t-1}}^{\textrm{bat}} + (\eta_{\textrm{bat}} \cdot  P_{\textrm{B,j,t}}^{\textrm{ch}} - \frac{P_{\textrm{B,j,t}}^{\textrm{dis}}}{\eta_{\textrm{bat}}}) \cdot \Delta t,    \label{eq:dynBESS1}\\
    P_{\textrm{B,j,t}}^{\textrm{ch}} +   P_{\textrm{B,j,t}}^{\textrm{dis}} \leq  \textrm{max}(P_{\textrm{B,j,t}}^{\textrm{ch}}, P_{\textrm{B,j,t}}^{\textrm{dis}}), \label{eq:BESS_P2} \\
    Q_{\textrm{B,j,t}}^{2} \leq  (S_{\textrm{max}}^{\textrm{bat}})^{2} - \textrm{max}((P_{\textrm{B,j,t}}^{\textrm{ch}})^{2}, (P_{\textrm{B,j,t}}^{\textrm{dis}})^{2}),
\end{gather}
\end{subequations}

\noindent where $E_{\textrm{j,t}}^{\textrm{bat}}$ is the available energy capacity at node $j$ and time $t$ with an initial \SK{and terminal} energy content $E_{\textrm{start}}$ and constrained by $SoC_{\textrm{min}}^{\textrm{bat}}$, $SoC_{\textrm{max}}^{\textrm{bat}}$ which are the fixed minimum and maximum per unit limits; $E_{\textrm{inv,j}}^{\textrm{bat}}$ is the installed BESS capacity at node $j$, $P_{\textrm{B,j,t}}^{\textrm{ch}}$ and $ P_{\textrm{B,j,t}}^{\textrm{dis}}$ are the charging and discharging powers defined as positive according to \eqref{eq:BESS_P}; the energy capacity at each time step $t$ influenced by the BESS efficiency $\eta_\textrm{bat}$ is defined by \eqref{eq:dynBESS1} accounting for the time interval $\Delta t$; \eqref{eq:BESS_P2} is re-casted as mixed-integer constraint with two binaries for each time step, and ensures that the BESS is not charging and discharging at the same time. In case this constraint leads to high computational burden, it can be substituted by
\begin{gather}
    P_{\textrm{B,j,t}}^{\textrm{ch}} \cdot  (P_{\textrm{l,j,t }} - P_{\textrm{g,j,t}}^{\textrm{max}} )  \leq \epsilon, \label{Ch3eq:charging} \\
    P_{\textrm{B,j,t}}^{\textrm{dis}} \cdot (P_{\textrm{l,j,t}} - P_{\textrm{g,j,t}}^{\textrm{max}} )  \geq \epsilon,  \label{Ch3eq:discharging}
\end{gather}
where we use an arbitrarily small value $\epsilon=10^{-5}$ to approximate the original constraint \eqref{eq:BESS_P2}, by not allowing the BESS to discharge (resp. charge) at times of excess (resp. deficit) local generation. \SK{This approximation is exact when using the current static dual tariff scheme in distribution grids, or when the consumer is aiming at maximizing the self-consumption~\cite{karagiannopoulos2019active}.}

\subsection{Inverter technical capability}
The increasing installation of DERs in MV and LV grids necessitates the modernization of the standards and grid codes. The main goal is to account for operation modes that exploit the DER controllability and flexibility capabilities to the benefit of both the system (TN or DN) and the customers. Although the initial requirements referred to MV DERs, it has become clear that the application to LV units will be very advantageous~\cite{karagiannopoulos2019data,StavrosPSCC18}.  
\begin{figure}[]
    \begin{centering}
	\includegraphics[width=0.99\columnwidth]{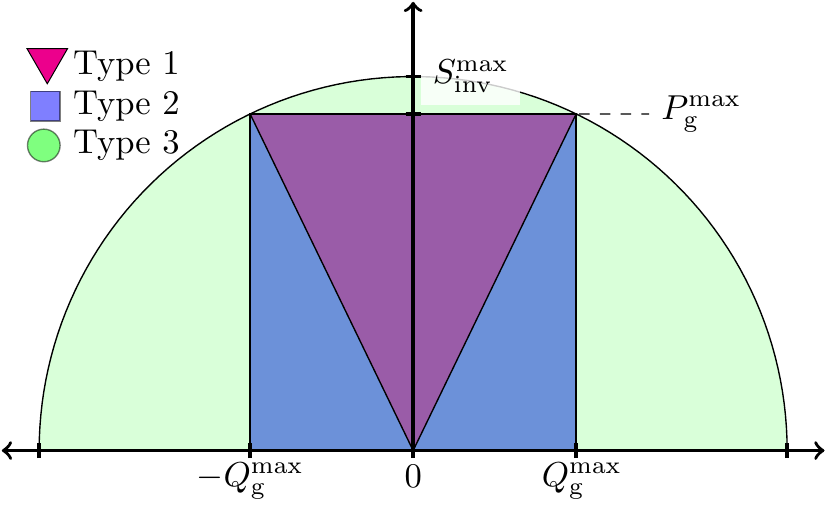}
	 \caption{Different operational types within the P-Q inverter capability curve.}	\label{fig:inverters} 
	\end{centering} 
\end{figure}
Regarding voltage support, several standards describe rules for DERs connecting in LV grids~\cite{VDE,ItalianNorm2014,IEEE1547}. Most of them prescribe operation with a minimum power factor which can be enforced by local control schemes using volt/var characteristic curves, a fixed power factor operation, or a fixed reactive power behavior. For example,~\cite{VDE} follows an open-loop local scheme in which the DERs adjust their power factor as a function of the active power they inject into the grid. The lowest power factor therein is reached when they inject power at their capacity limit and the values are $\cos\phi=0.95$ and $\cos\phi=0.9$ for small and for large PV capacities, respectively. 
In the USA, only the recently revised version of IEEE 1547~\cite{IEEE1547} allows DERs on primary and secondary distribution voltages to provide steady-state voltage support. Closed-loop schemes based on volt/var curves result typically in a more efficient grid utilization, since they link the additional reactive power consumption with instances of high local voltages. However, since they are allowed only in a few standards, e.g.~\cite{ItalianNorm2014,VDE,Bdew2008}, the need for a harmonization of the allowed operation modes emerges. The interested reader is referred to~\cite{Kotsampopoulos2013} for a comparison of existing guidelines in terms of DER control. Our earlier work~\cite{karagiannopoulos2019data} focused on providing optimal local curves instead of such standardized curves. 
However, if a reliable monitoring and communication infrastructure is available, optimal coordination is feasible in close to real-time without the need for designing local decision curves.

Inverter-based DERs can operate within the P-Q inverter capability curve, depicted in Fig.~\ref{fig:inverters}. In many cases, the inverter is overdimensioned by $5\%$ or $10\%$, e.g. $S_{\textrm{inv,j}}^{\textrm{max}}=1.1 \cdot P_{\textrm{g,j}}^{\textrm{max}}$ in order to allow for reactive power control even when the DER is operating at the maximum active power. \SK{In this paper, we will examine the following cases: a) the ``triangular" ($\bigtriangledown$) limitations of~\eqref{eq:PQ_trigwno} that are used in~\cite{Bdew2008,VDE,ItalianNorm2014}, and impose an operational minimum power factor linking the active and reactive power injections, b) the ``rectangular" ($\square$) limitations of~\eqref{eq:PQ_square} defined in~\cite{ItalianNorm2014} which allows for reactive power control even at times with low active power injections, and finally c) the semi-circle ($\bigcirc$) capability described by~\eqref{eq:PQ_semicirc} which illustrates the full capability region of the DERs, without imposing any constraint on the operational power factor.} The latter is possible because lately, there are schemes that offer reactive power capabilities during inactive times, e.g. at night when the PV units do not inject power~\cite{SMA}. The equations for the considered cases are given by

\begin{subequations}
\label{eq:PQlimits}
\begin{align}
	(\bigtriangledown):& -\textrm{tan}(\phi_{\textrm{max}})  P_{\textrm{g,j,t}}^{\textrm{ }}  \leq  Q_{\textrm{g,j,t}}^{\textrm{ }} \leq \textrm{tan}(\phi_{\textrm{max}})  P_{\textrm{g,j,t}}^{\textrm{ }}, \label{eq:PQ_trigwno} \\
 	(\square):& -\textrm{tan}(\phi_{\textrm{max}})  P_{\textrm{g,j,t}}^{\textrm{min}} \leq  Q_{\textrm{g,j,t}}^{\textrm{ }} \leq \textrm{tan}(\phi_{\textrm{max}})  P_{\textrm{g,j,t}}^{\textrm{max}}, \label{eq:PQ_square} \\
 	(\bigcirc):& \qquad \qquad Q_{\textrm{g,j,t}}^{2} \leq  (S_{\textrm{inv,j}}^{\textrm{max}})^{2} - P_{\textrm{g,j,t}}^{2}. \label{eq:PQ_semicirc} 
\end{align}
\end{subequations}

\section{Case Study - Results}\label{case}

\begin{figure}
    \begin{centering}
	\includegraphics[width=1\columnwidth]{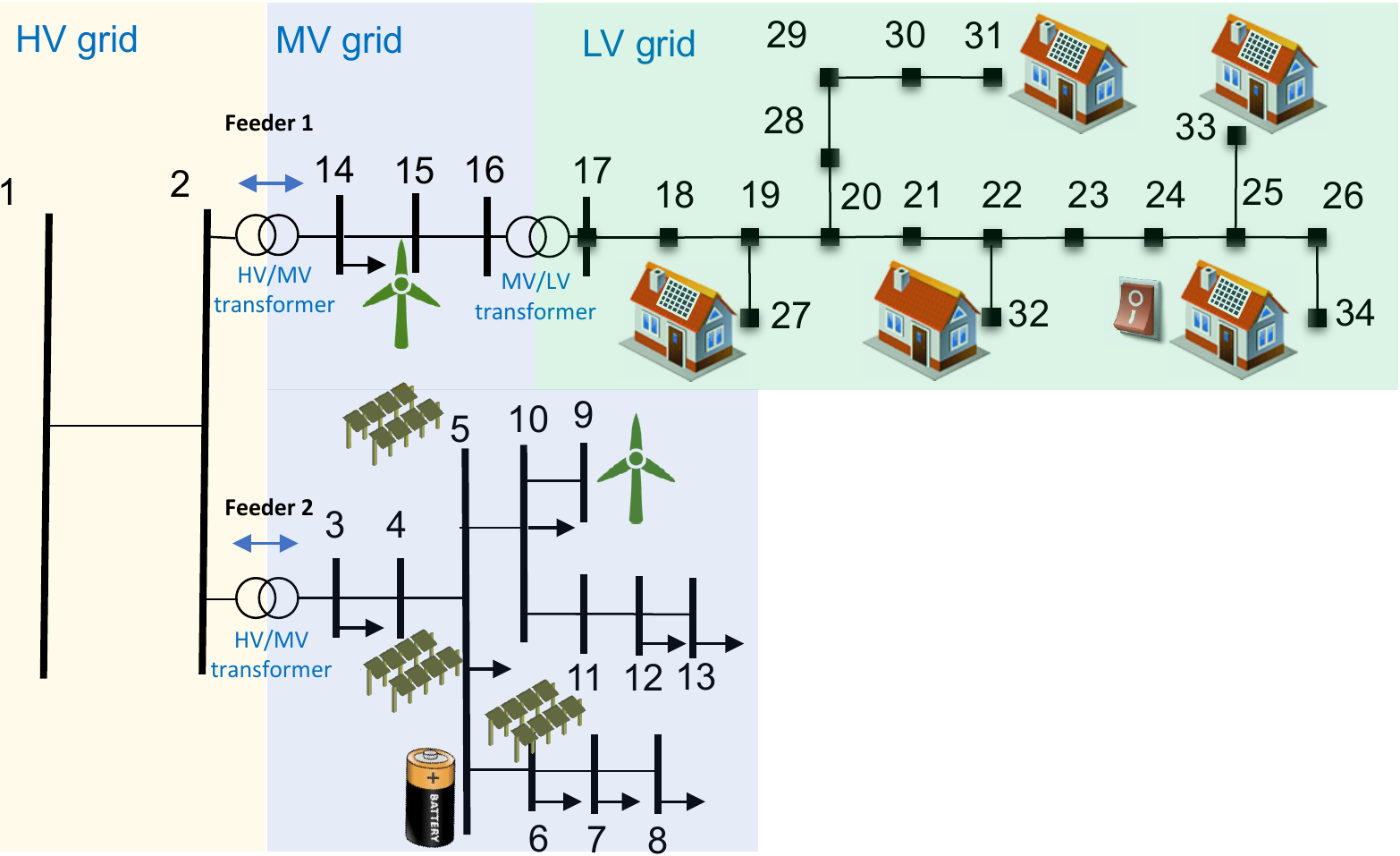}
	 \caption{Joint MV and LV European grids~\cite{Strunz2014} connected to the transmission voltage level through a Th\'{e}venin equivalent.}	\label{fig:cigre_test_system} 
	\end{centering}
\end{figure}

To analyze the performance of the proposed centralized scheme, we merged the typical European MV and LV grids~\cite{Strunz2014}, as sketched in Fig.~\ref{fig:cigre_test_system}. More specifically, in contrast to~\cite{Strunz2014} which models each voltage level separately, we constructed the MV benchmark and substituted one of the loads with the LV residential grid. Thus, we consider two feeders connecting to the HV network. In the first feeder, we model in detail the LV grid with the operational flexibility of each rooftop PV unit, while the second aggregates the LV grid similar to most studies. The TN is modeled by a Th\'{e}venin equivalent at Node $1$ in order to investigate the influence of the distribution grid on the voltage support scheme. \SK{The parameters of the Th\'{e}venin equivalent used  in this work are very important and define the effectiveness of voltage support from the distribution grid. In stiff HV grids, high amounts of reactive power are needed to control voltages while in weak grids low levels of reactive power already have impact. Due to the energy deregulation requirement, DSOs cannot have parts of the actual model of the transmission networks or neighboring distribution grids which might belong to other operators. In fact, in many countries the TSO calculates and communicates to the DSOs the Th\'{e}venin impedance and the voltage reference values at each interconnection point. Having forecasts for the load and the generation, the TSO can estimate these values based on the hourly load and generation mix. In Switzerland, these values do not change significantly when compared to weak grids, e.g. of an island. Thus, a relatively high value for the Swiss system is selected to represent a weak future 220 kV grid with low inertia, i.e. three times the inverse of the short circuit capacity (SCC) available on the primary side of the transformer according to~\cite{Strunz2014}. Different values can be incorporated easily in our formulation, since it corresponds to different hourly values.}

The installed PV and wind turbine capacities in the MV and LV grids are summarized as follows: PV nodes~$=[4, 5, 6, 27, 31, 33, 34]$, PV~installed~capacity (MVA)~$=[1, 1, 1, 0.1, 0.15, 0.1, 0.12]$, WT nodes~$=[9, 15]$, WT~installed~capacity (MVA)~$=[1, 1]$. \SK{Furthermore, we consider at node~$5$ a BESS of $200$~kWh, and at node $34$ a CL of $5$ kW, whose total daily energy consumption needs to be maintained constant, to highlight the potential coordinated contribution of all available DERs.} 
In this work, we only consider balanced, single-phase system operation, but the framework can be extended to three-phase unbalanced networks as in~\cite{StavrosPSCC18,karagiannopoulos2019data}. The normalized PV injection profiles are taken from PV stations in Switzerland, following~\cite{StavrosIREP}, the wind profiles from~\cite{SmartPlanningWP2.5} and the load data are taken from~\cite{Strunz2014}. The assumed installed capacities result in overvoltage issues when all DGs are injecting power as is the case in reality with high PV penetrations.  

The operational costs are set to $C_{\textrm{curt}} = 0.3 \frac{\textrm{CHF}}{\textrm{kWh}}$ and $C_{\textrm{Q}} = 0.01 \cdot C_{\textrm{curt}}$. The BESS cost is considered in the planning stage~\cite{stavrosPowertech} and thus, the use of the BESS does not incur any operational cost to the DNO. The tariffs for the active and passive participation in voltage support for 2018 are set to $C_{\textrm{p}} = 0.0151 \frac{CHF}{kvarh}$, $c_{\textrm{a,c}} = 0.003 \frac{CHF}{kvarh}$ and $c_{\textrm{a,n}} = 0.0151 \frac{CHF}{kvarh}$ as given in~\cite{SwissGrid_VoltageSupport_Tariff}.

The implementation was done in MATLAB. For the centralized OPF-based control, YALMIP~\cite{Lofberg2004} was used as the modeling layer and Gurobi~\cite{gurobi} as the solver. The results were obtained on an Intel Core i7-2600 CPU and 16 GB of RAM.

\subsection{Passive participation}\label{pass_parti_case}
\begin{figure}[t]
    \begin{centering}
	\includegraphics[width=0.99\columnwidth]{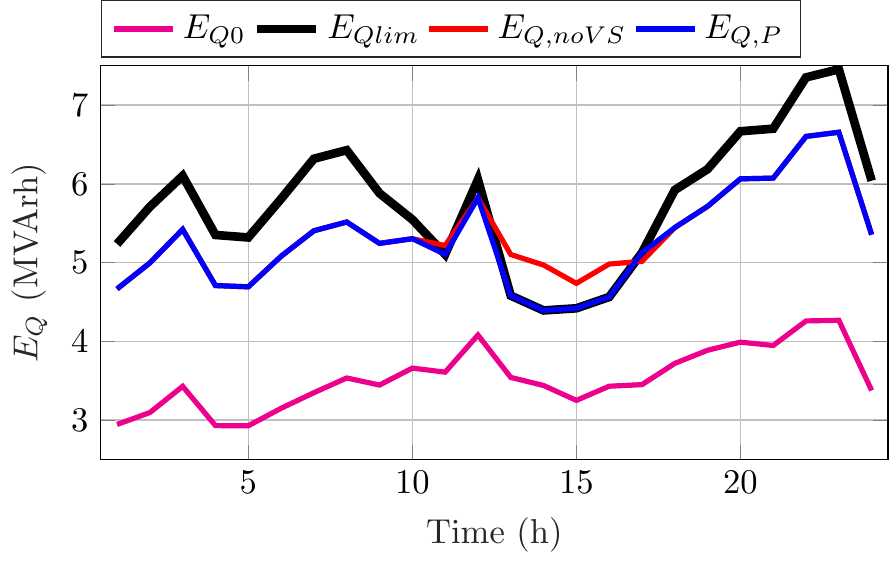}
	 \caption{Daily reactive power exchange with the TN of Feeder~1.}	\label{fig:LowTr} 
	\end{centering} 
\end{figure}
\SK{ In the first part of the results we examine the case, where the DN follows a passive participation in the voltage support scheme. We provide detailed results for a specific day, i.e. $96$ timesteps, and then summarized seasonal results. \subsubsection{Daily Results} 
We compare} the results from a) the no-control case where all DERs inject the maximum available active power without reactive power control, b) running an OPF-based scheme as described in Section II without providing voltage support to the TN, i.e. excluding the constraints referring to the compliant region and the respective objective function term, and c) running the same OPF-based scheme, following in addition the passive voltage support participation. Figure~\ref{fig:LowTr} shows the daily reactive energy exchange between Feeder $1$ of the DN and the TN. First, we observe that the daily reactive energy exchange without controlling any DER ($E_{\textrm{Q0}}$), i.e. just computing the power flows, is lower than the profiles of the OPF but without providing voltage support ($E_{\textrm{Q,noVS}}$) and the OPF including passive voltage support participation ($E_{\textrm{Q,P}}$). That means, that the OPF manages to satisfy all network constraints, by requiring more reactive power from the TN to impose an inductive behavior from some DERs and hence, eliminate overvoltage issues. 

\begin{figure}[b]
    \begin{centering}
	\includegraphics[width=0.99\columnwidth]{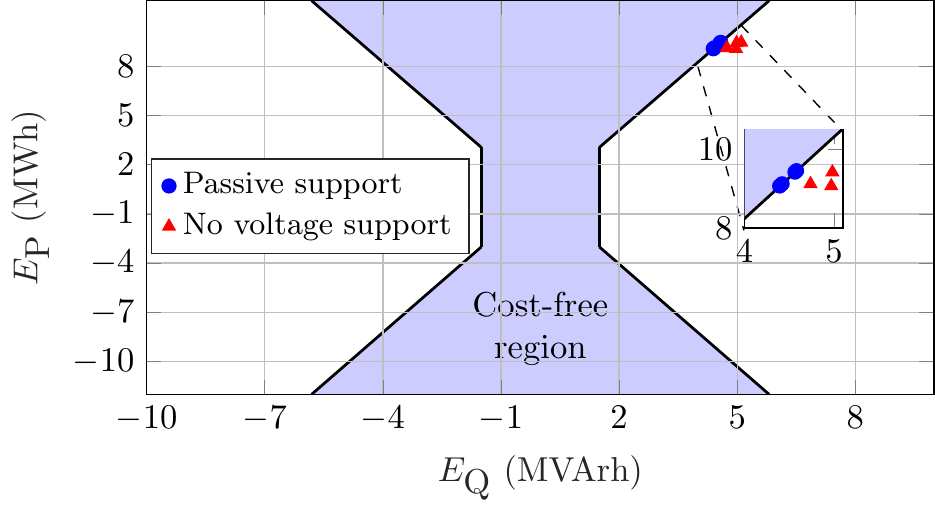}
	 \caption{Cost-free and penalization regions for the passive voltage support participation.}	\label{fig:passive_case} \vspace{-0.3cm}
	\end{centering} 
\end{figure}

The voltage support scheme incurs costs when the reactive power exchange exceeds the varying limit ($E_{\textrm{Q,noVS}}$). As can be seen, during noon hours the OPF that considers the voltage support scheme manages to reach the cost-free region, by requiring less reactive power from the TN. This does not result in any constraint violation, since other measures are activated to satisfy the operational constraints. The change of the operating points for the hours 13:00, 14:00 and 15:00 are shown in Fig.~\ref{fig:passive_case}. The DN is able to reduce only the reactive energy exchange in order to avoid the cost penalization without interfering with the active energy exchange.    

\SK{\subsubsection{Seasonal economic evaluation} In this part, we compare the total cost over the period of 3 months. Table~\ref{tab:costPassive} summarizes the comparison of the total costs for the summer months which show high injections from the PV parks and panels. Generally, the cost for the voltage support scheme is much lower than the other objective function terms, e.g. curtailment and losses. Thus, there are days where the voltage support scheme does not show as large cost savings as in the day presented in the daily analysis. However, over the time period of 3 summer months, the consideration of the passive voltage support scheme in the OPF formulation resulted in a reduction of around $3 kCHF$, corresponding in a $3.39\%$ cost reduction compared to the case with no voltage support. This value depends on the operator's tariffs and the trade-offs among all the costs in the objective function.}

\begin{table}[] 
\caption{Seasonal cost comparison for the passive role.}
\centering
\begin{tabular}{c|c|c}
\toprule
Type of Control & \begin{tabular}[c]{@{}c@{}}Total Seasonal\\ Cost (kCHF)\end{tabular} & \begin{tabular}[c]{@{}c@{}}Cost Difference \\ due to voltage\\ support (\%)\end{tabular}\\ \midrule
\begin{tabular}[c]{@{}c@{}}OPF-based: \\ No voltage support\end{tabular}         & $98.78$                      & - (reference)        \\ \hline
\begin{tabular}[c]{@{}c@{}}OPF-based:\\ Passive role\end{tabular} & $95.42$                      & $-3.39$                \\ \bottomrule
\end{tabular}
\label{tab:costPassive}
\end{table}

\subsection{Active participation}\label{active_parti_case}
\subsubsection{Daily Results}
In this part, the DN participates in active voltage support as explained in Section~\ref{activeVscheme}. The focus here is to use the operational flexibility of the DERs located in both MV and LV in order to track a voltage magnitude profile provided by the TSO. Furthermore, we highlight the importance of an extended operational range for reactive power control by providing results for all three types of the inverter active-reactive power capability curve. That is a) the ``triangular" ($\bigtriangledown$) limitations of~\eqref{eq:PQ_trigwno}, the ``rectangular" ($\square$) limitations of~\eqref{eq:PQ_square} and finally c) the semi-circle ($\bigcirc$) limitations of~\eqref{eq:PQ_semicirc}. The voltage reference profile which is derived by the TSO DARP procedure, is assumed flat at $1$ p.u. with a $0.5\%$ tolerance, i.e. $\epsilon=0.005$ p.u..  

\begin{figure}[]
    \begin{centering}
	\includegraphics[width=0.99\columnwidth]{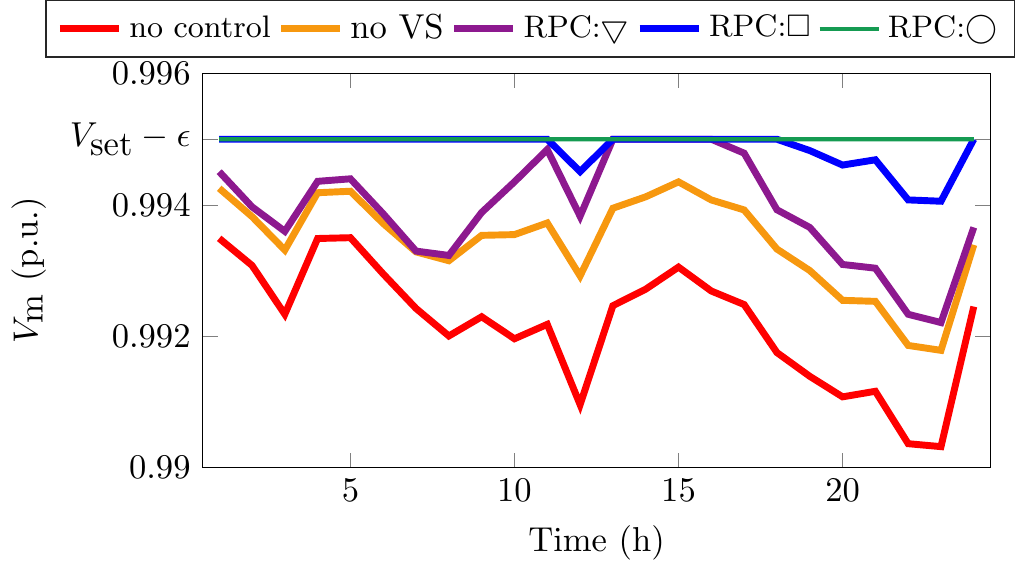}
	 \caption{Daily voltage profile for the centralized schemes participating in the active voltage support scheme.}	\label{fig:active_voltages} \vspace{-0.3cm}
	\end{centering} 
\end{figure}

Figure~\ref{fig:active_voltages} shows the daily voltage magnitude profile for the different cases. Apart from the OPF results we first provide the ``no control" case, by running simple power flow calculations without any control, i.e. power factor set to one for all PVs. We observe that the OPF case that does not consider the active voltage support scheme, i.e. ``no VS case",  results in a voltage increase, indicating that less reactive power is required by the TN and more reactive power is produced locally within the DN, according to the objective function minimization of ~\eqref{eq:objfun1}. However, the voltage needs to be further increased in order to reach the compliant area defined by the active voltage support scheme. The ``triangular" case ($\bigtriangledown$) shows limited capability, being linked with the active power injection of the DERs. Over the night hours, the PV units are not active and it is only the flexibility of the WTs that can contribute to the voltage raise. However, the ``rectangular" ($\square$) case manages to reach the complaint region for most hours of the day, but only the semi-circle variant is remunerated for the whole day since it exploits the largest inverter capability. 

\begin{figure}[]
    \begin{centering}
	\includegraphics[width=0.99\columnwidth]{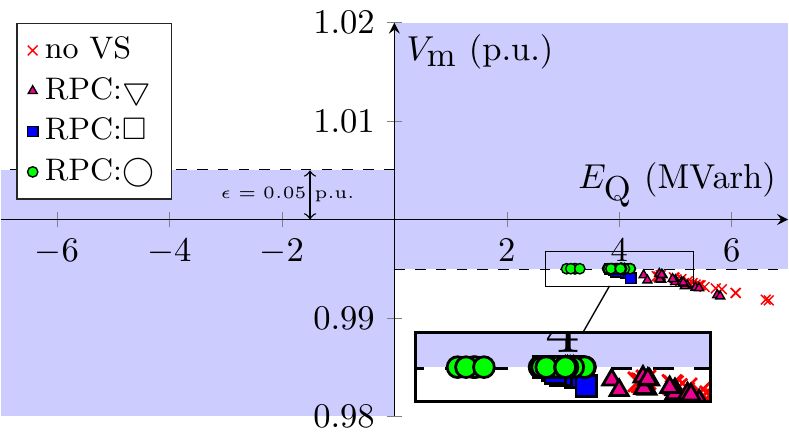}
	 \caption{Compliant and penalization regions for the active voltage support participation.}	\label{fig:active_caseRegions} \vspace{-0.3cm}
	\end{centering} 
\end{figure}

Figure~\ref{fig:active_caseRegions} shows the pairs of the measured voltage and the daily reactive energy exchange for the different cases with a time resolution of $1$ hour. As we observe, the larger the reactive power capabilities of the DER inverters, the closer the DN can get to the compliant regions. 

\begin{figure}[]
    \begin{centering}
	\includegraphics[width=0.99\columnwidth]{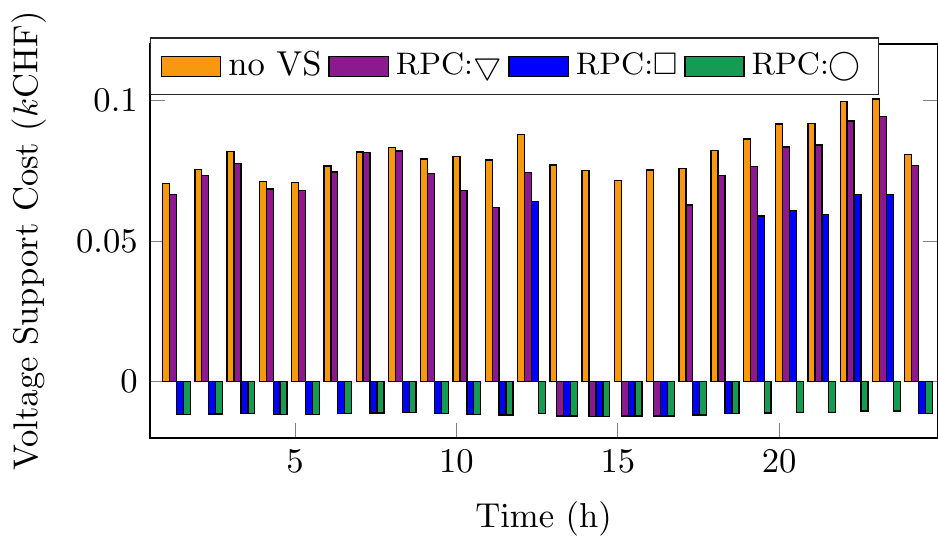}
	 \caption{Hourly cost and revenue for the active voltage support participation cases.}	\label{fig:costbars} \vspace{-0.3cm}
	\end{centering} 
\end{figure}

Finally, Fig.~\ref{fig:costbars} shows the hourly cost and revenues for the different variations of the active voltage support participation. As it is observed, the ``no VS" case results in the highest costs for each hour. Then, the larger the P-Q capability region, the lower the costs, with the semi-circle ($\bigcirc$) showing small revenues at all hours of the day.

It is important to note that simulating the LV grid guarantees that the power quality is not compromised in the LV level due to actions in other voltage levels. This allows us to quantify the real potentials of each unit according to its location and technical capabilities, and finally unlocks business cases for aggregators in distribution grid levels. 


\SK{Seasonal economic evaluation: Similar to the passive role, here we compare the total cost over the period of 3 months following the active participation role. Table~\ref{tab:costActive} summarizes the voltage support costs for the summer months considering the three alternatives in terms of inverter capabilities. First, we observe that these costs are much higher than the ones of the passive case, since there is no cost-free region anymore. The costs decrease with higher operational flexibility and, similar to the daily results, only the semi-circle case results in revenues. The rectangular capability succeeds in reducing the voltage support cost by around $67\%$ more than the triangular case, highlighting the benefits of controlling reactive power even when active power is very low.}
\begin{table}[t]
\caption{Seasonal voltage support cost for the active role.}
\centering
\begin{tabular}{c|c|c}
\toprule
\begin{tabular}[c]{@{}c@{}}Type of Control\\ (OPF-based)\end{tabular} & \begin{tabular}[c]{@{}c@{}}Seasonal cost \\ for Voltage\\ support (kCHF)\end{tabular} & \begin{tabular}[c]{@{}c@{}}Cost Difference \\ due to voltage\\ support (\%)\end{tabular} \\ \midrule
No voltage support                                                    & 52.20                                                                                 & (reference)                                                                            \\ \hline
Active role ($\bigtriangledown$)                                                         & 39.56                                                                                 & -24.65                                                                                   \\ \hline
Active role ($\square$)                                                         & 4.51                                                                                  & -91.41                                                                                   \\ \hline
Active role ($\bigcirc$)                                                        & -7.42                                                                                 & -114.14                                                                                  \\ \bottomrule
\end{tabular}
\label{tab:costActive}
\end{table}

\section{Conclusion}  \label{Conclusion}
The increasing controllability and observability in MV and LV grids call upon a more efficient DN-TN collaboration. Active DNs are capable of providing ancillary services to the transmission voltage levels, upgrading their conventional role from power sinks to flexible participants. Controlling actively different types of DERs enhances the security, resilience and optimal operation of multiple voltage levels. In this paper, we propose a tractable centralized OPF methodology to provide voltage support to the TN, while optimizing the DN operation and satisfying the power quality constraints. We have demonstrated through case studies that in the presence of communication and monitoring infrastructure, a DSO can optimize its grid operation safely, while at the same time profit from voltage support schemes that yield revenues when contributing to the TSO objectives. \SK{In the considered systems, the new active role participation increases the voltage support costs significantly compared to the passive case, motivating the DSOs to utilize their reactive power flexibility. In the active role, the rectangular operational type of the P-Q inverter capability led to a significant voltage cost reduction. Finally}, we have highlighted the importance of exploiting the full reactive power potential of the inverter-based DERs, pointing out the need for modernization and harmonization of grid codes and standards to allow more advanced DER control in DNs. 
\bibliographystyle{IEEEtran}
\bibliography{bibliography}

\vspace{-0.13cm}

\end{document}